\documentclass[12pt,oneside,english]{amsart}
\usepackage[T1]{fontenc}
\usepackage[latin9]{inputenc}
\usepackage{geometry}
\usepackage{color}
\usepackage{babel}
\usepackage{verbatim}
\usepackage{amsthm}
\usepackage{amssymb}
\usepackage{esint}
\usepackage[unicode=true]
 {hyperref}
\usepackage{breakurl}

\makeatletter
\numberwithin{equation}{section}
\numberwithin{figure}{section}
\theoremstyle{plain}
\newtheorem{thm}{\protect\theoremname}[section]
  \theoremstyle{remark}
  \newtheorem{rem}[thm]{\protect\remarkname}
  \theoremstyle{plain}
  \newtheorem{lem}[thm]{\protect\lemmaname}
  \theoremstyle{plain}
  \newtheorem{fact}[thm]{\protect\factname}

\makeatother

  \providecommand{\factname}{Fact}
  \providecommand{\lemmaname}{Lemma}
  \providecommand{\remarkname}{Remark}
\providecommand{\theoremname}{Theorem}

\begin{document}

\title[The Peetre ''plus-minus'' interpolation spaces $\left\langle A_{0},A_{1}\right\rangle _{\theta}$]{Some alternative definitions for the ''plus-minus'' interpolation
spaces $\left\langle A_{0},A_{1}\right\rangle _{\theta}$ of Jaak
Peetre}

\author{Michael Cwikel}

\address{Department of Mathematics, Technion - Israel Institute of Technology,
Haifa 32000, Israel }

\email{mcwikel@math.technion.ac.il }

\thanks{Research supported by the Technion V.P.R.\ Fund and by the Fund
for Promotion of Research at the Technion.}
\begin{abstract}
The Peetre ``plus-minus'' interpolation spaces $\left\langle A_{0},A_{1}\right\rangle _{\theta}$
are defined via conditions regarding the unconditional convergence
of Banach space valued series of the form $\sum_{n=-\infty}^{\infty}2^{(j-\theta)n}a_{n}$
or, alternatively, $\sum_{n=-\infty}^{\infty}e^{(j-\theta)n}a_{n}$,
for $j\in\left\{ 0,1\right\} $. It may seem intuitively obvious that
using powers of $2$ or of $e$, or powers of some other constant
number greater than $1$ in these definitions should produce the same
spaces to within equivalence of norms. To allay any doubts, we here
offer an explicit proof of this fact, via a ``continuous'' definition
of the same spaces, where integrals replace the above series. This
apparently new definition, which is also in some sense a ``limiting
case'' of the usual ``discrete'' definitions, may be relevant in
the study of the connection between $\left\langle A_{0},A_{1}\right\rangle _{\theta}$
and the Calder\'on complex interpolation space $\left[A_{0},A_{1}\right]_{\theta}$
in the case where $\left(A_{0},A_{1}\right)$ is a couple of Banach
lattices. Related results can probably be obtained for the Gustavsson-Peetre
variant of the ``plus-minus'' spaces.
\end{abstract}
\maketitle

\section{\label{sec:Intro}``Pre-Introduction'' - some background}

In \cite{PeetreJ1971}, among a number of other very interesting things,
Jaak Peetre introduced several kinds of interpolation spaces, including
ones which we will soon define in Section \ref{sec:BasicProps}, and
which he denoted by $\left\langle A_{0},A_{1}\right\rangle _{\theta}$.
(The bibliography of \cite{PeetreJ1971} indicates that such spaces,
or some variants of them, had apparently already been considered as
early as in 1962 and were discussed in (Swedish) lecture notes written
jointly at that time by Peetre and Arne Persson, for a course at the
University of Lund.) 

It should be mentioned that, a decade after the appearance of \cite{PeetreJ1971},
the remarkable paper \cite{JansonS1981} of Svante Janson identified
$\left\langle A_{0},A_{1}\right\rangle _{\theta}$ as an ``orbit''
space and offered other additional insights about it (and about several
other kinds of interpolation spaces). Continuing in the spirit of
\cite{PeetreJ1971}, the two papers \cite{JansonS1981} and \cite{CwikelMKaltonNkmr}
each exhibited different ways of fitting the spaces $\left\langle A_{0},A_{1}\right\rangle _{\theta}$
and other kinds of interpolation spaces into more general frameworks.
(The approach of \cite{JansonS1981} is elaborated upon further in
\cite{OvchinnikovV1984}, and some further discussion of the contents
of \cite{CwikelMKaltonNkmr} can be found in \cite{CwikelMMilmanMRochbergR2014B}.)
But the only result from these papers which I will use here is a counterexample
from \cite{JansonS1981} which will be mentioned below in Remark \ref{rem:SvanteExample}.
Also I will not deal with the variants $\left\langle A_{0},A_{1},\rho\right\rangle $
of the spaces $\left\langle A_{0},A_{1}\right\rangle _{\theta}$ which
were introduced and studied in \cite{GustavssonJPeetreJ1977} by Jan
Gustavsson and Peetre and also later in \cite{JansonS1981}. The contents
of this note can surely be adapted to give similar results about the
spaces $\left\langle A_{0},A_{1},\rho\right\rangle $, in particular
when the function parameter $\rho$ in their definition is of the
form $\rho(t)=t^{\theta}$ for some $\theta\in(0,1)$.

The contents of the rest of this note are as follows: In Section \ref{sec:BasicProps}
we recall the various slightly different ``discrete'' definitions
of the spaces $\left\langle A_{0},A_{1}\right\rangle _{\theta}$,
i.e., the definitions which have been used to date, and which have
been implicitly assumed to be equivalent to each other. We also mention
some of the basic properties of these spaces and of the auxiliary
spaces which are used to define them. In Section \ref{sec:CtsDefn}
we introduce a new ``continuous'' definition of $\left\langle A_{0},A_{1}\right\rangle _{\theta}$
and show that it is equivalent to each of the ``discrete'' definitions
of Section \ref{sec:BasicProps}, and therefore we also confirm that
those definitions are indeed equivalent to each other. Section \ref{sec:PeetreInCalderon}
could be considered to be a sort of appendix. In it we recall the
well known and obvious fact that the space $\left\langle A_{0},A_{1}\right\rangle _{\theta}$
is contained in the complex interpolation space $\left[A_{0},A_{1}\right]_{\theta}$,
and note that this is a norm one inclusion. Section \ref{sec:Appendix}
is most definitely an appendix, which describes straightforward proofs
of the completeness of $\left\langle A_{0},A_{1}\right\rangle _{\theta}$
and of the auxiliary space used to define it, for the convenience
of any reader who is less familiar with these matters.

\section{\label{sec:BasicProps} More introduction - the standard ``discrete''
definition(s) of the space $\left\langle A_{0},A_{1}\right\rangle _{\theta}$,
and some of its properties}

We recall that a two-sided series $\sum_{n\in\mathbb{Z}}a_{n}$ of
elements $a_{n}$ in a Banach space $A$ is said to be \textit{unconditionally
convergent} if, for every bounded sequence $\left\{ \lambda_{n}\right\} _{n\in\mathbb{Z}}$
of scalars, the sequence of partial sums $\left\{ \sum_{n=-N}^{N}\lambda_{n}a_{n}\right\} _{N\in\mathbb{N}}$
converges in norm in $A$ as $N$ tends to $\infty$. This property
readily implies, and is therefore equivalent to, the following stronger
condition: 
\begin{equation}
\lim_{N\to\infty}\sup\left\{ \left\Vert \sum_{\left|n\right|\ge N}\lambda_{n}a_{n}\right\Vert _{A}:\left\{ \lambda_{k}\right\} _{k\in\mathbb{Z}}\in\Lambda\right\} =0,\label{eq:LikeCzero}
\end{equation}
\label{pCzero}where $\Lambda$ is the set of all scalar sequences
$\left\{ \lambda_{k}\right\} _{k\in\mathbb{Z}}$ which satisfy $\sup_{k\in\mathbb{Z}}\left|\lambda_{k}\right|\le1$.
(The straightforward proof that unconditional convergence implies
(\ref{eq:LikeCzero}) is a simpler variant of the proof below of Fact
\ref{fact:zPrinch}.) 

Of course (\ref{eq:LikeCzero}) implies that 
\begin{equation}
\sup\left\{ \left\Vert \sum_{n=-\infty}^{\infty}\lambda_{n}a_{n}\right\Vert _{A}:\left\{ \lambda_{k}\right\} _{k\in\mathbb{Z}}\in\Lambda\right\} <\infty.\label{eq:bpi}
\end{equation}

\begin{rem}
Note that, in seeming contrast to some phenomena in some other more
elementary situations, the existence of the limits $\lim_{N\to\infty}\sum_{n=-N}^{N}\lambda_{n}a_{n}$
for all $\left\{ \lambda_{k}\right\} _{k\in\mathbb{Z}}\in\Lambda$
is obviously equivalent to the existence of both of the limits $\lim_{N\to\infty}\sum_{n=1}^{N}\lambda_{n}a_{n}$
and $\lim_{N\to\infty}\sum_{n=-N}^{0}\lambda_{n}a_{n}$ for all $\left\{ \lambda_{k}\right\} _{k\in\mathbb{Z}}\in\Lambda$.
\end{rem}
Let $\mathbb{F}$ be either $\mathbb{R}$ or $\mathbb{C}$ and let
$\left(A_{0},A_{1}\right)$ be an arbitrary Banach couple where $\mathbb{F}$
is the field of scalars for $A_{0}$ and also for $A_{1}$. For each
choice of the constants $r>0$ and $\theta\in(0,1)$, let $J(\theta,r,A_{0},A_{1})$
denote the space of $A_{0}\cap A_{1}$ valued sequences $\left\{ a_{n}\right\} _{n\in\mathbb{\mathbb{Z}}}$
having the two properties that $\sum_{n\in\mathbb{Z}}e^{-r\theta n}a_{n}$
is an unconditionally convergent series in $A_{0}$ and $\sum_{n\in\mathbb{Z}}e^{r(1-\theta)n}a_{n}$
is an unconditionally convergent series in $A_{1}$. In view of the
properties of unconditionally convergent series discussed above, we
see (cf.~(\ref{eq:bpi})) that the quantity 
\begin{equation}
\left\Vert \left\{ a_{n}\right\} _{n\in\mathbb{Z}}\right\Vert _{J(\theta,r,A_{0},A_{1})}:=\sup\left\{ \left\Vert \sum_{n=-\infty}^{\infty}\lambda_{n}e^{(j-\theta)rn}a_{n}\right\Vert _{A_{j}}:j\in\left\{ 0,1\right\} ,\,\left\{ \lambda_{k}\right\} _{k\in\mathbb{Z}}\in\Lambda\right\} \label{eq:defjn}
\end{equation}
is finite for every sequence $\left\{ a_{n}\right\} _{n\in\mathbb{Z}}$
in $J(\theta,r,A_{0},A_{1})$ and therefore defines a norm on that
space.
\begin{rem}
\label{rem:LambdaZero}It is clear that, for each $\left\{ a_{n}\right\} _{n\in\mathbb{Z}}$
in $J(\theta,r,A_{0},A_{1})$, the value of the supremum in (\ref{eq:defjn})
will not change if we replace $\Lambda$ in (\ref{eq:defjn}) by its
subset $\Lambda_{0}$ consisting of those sequences which have only
finitely many non-zero elements.
\end{rem}

\begin{rem}
\label{rem:DiscreteMonotone}We note that, for $j\in\left\{ 0,1\right\} $
and every sequence $\left\{ a_{n}\right\} _{n\in\mathbb{Z}}$ in $J(\theta,r,A_{0},A_{1})$,
\begin{equation}
\sup\left\{ \left\Vert \sum_{n\in U}\lambda_{n}e^{(j-\theta)rn}a_{n}\right\Vert _{A_{j}}:\left\{ \lambda_{k}\right\} _{k\in\mathbb{Z}}\in\Lambda\right\} \le\sup\left\{ \left\Vert \sum_{n\in V}\lambda_{n}e^{(j-\theta)rn}a_{n}\right\Vert _{A_{j}}:\left\{ \lambda_{k}\right\} _{k\in\mathbb{Z}}\in\Lambda\right\} \label{eq:itijgts}
\end{equation}
whenever $U$ and $V$ are (finite or infinite) sets which satisfy
$U\subset V\subset\mathbb{Z}$. This follows from the simple fact
that the left side of (\ref{eq:itijgts}) equals 
\[
\sup\left\{ \left\Vert \sum_{n\in V}\lambda_{n}e^{(j-\theta)rn}a_{n}\right\Vert _{A_{j}}:\left\{ \lambda_{k}\right\} _{k\in\mathbb{Z}}\in\Lambda_{U}\right\} ,
\]
where $\Lambda_{U}$ is the subset of $\Lambda$ consisting of all
those sequences $\left\{ \lambda_{k}\right\} _{k\in\mathbb{Z}}$ which
satisfy $\lambda_{k}=0$ for all $k\in\mathbb{Z}\setminus U$. 
\end{rem}

\begin{rem}
\label{rem:CauchyDef}For our convenience in some later parts of this
discussion, let us here explicitly state that, for every $A_{0}\cap A_{1}$
valued sequence $\left\{ a_{n}\right\} _{n\in\mathbb{N}}$, the following
three conditions are equivalent:

$\mathrm{(i)}$~ $\left\{ a_{n}\right\} _{n\in\mathbb{Z}}$ is an
element of $J(\theta,r,A_{0},A_{1})$.

$\mathrm{(ii)}$~ For $j\in\left\{ 0,1\right\} $ and for each $\varepsilon>0$
and each $\left\{ \lambda_{k}\right\} _{k\in\mathbb{Z}}\in\Lambda$,
there exists a positive integer $N(\varepsilon)$ (which could in
principle also depend on $\left\{ \lambda_{k}\right\} _{k\in\mathbb{Z}}$
and $j$) such that $\left\Vert \sum_{n_{1}\le\left|n\right|\le n_{2}}\lambda_{n}e^{(j-\theta)rn}a_{n}\right\Vert _{A_{j}}\le\varepsilon$
whenever $N(\varepsilon)\le n_{1}<n_{2}$. 

$\mathrm{(iii)}$~ For each $\varepsilon>0$, there exists a positive
integer $N(\varepsilon)$ such that 
\[
\left\Vert \sum_{n_{1}\le\left|n\right|\le n_{2}}\lambda_{n}e^{(j-\theta)rn}a_{n}\right\Vert _{A_{j}}\le\varepsilon
\]
 for every $\left\{ \lambda_{k}\right\} _{k\in\mathbb{Z}}\in\Lambda$
and for $j\in\left\{ 0,1\right\} $ whenever $N(\varepsilon)\le n_{1}<n_{2}$. 

The equivalence of (i) and (ii) and the implication (iii)$\Rightarrow$(ii)
are obvious. The implication (i)$\Rightarrow$(iii) follows immediately
from the remarks at the beginning of this section concerning the condition
(\ref{eq:LikeCzero}), i.e., from a straightforward proof similar
to that of Fact \ref{fact:zPrinch}.
\end{rem}
It is very easy to check (via reasoning similar to that which appears
below in Fact \ref{fact:Aplus}) that, for each $\left\{ a_{n}\right\} _{n\in\mathbb{Z}}\in J(\theta,r,A_{0},A_{1})$,
and each constant integer $k_{0}$, the partial sums $\sum_{n=-N}^{N+k_{0}}a_{n}$
converge in the norm of $A_{0}+A_{1}$, to a limit which does not
depend on $k_{0}$ and which we will denote by $\sum_{n=-\infty}^{\infty}a_{n}$,
and, furthermore, that this limit satisfies 
\begin{equation}
\left\Vert \sum_{n=-\infty}^{\infty}a_{n}\right\Vert _{A_{0}+A_{1}}\le2\left\Vert \left\{ a_{n}\right\} _{n\in\mathbb{Z}}\right\Vert _{J(\theta,r,A_{0},A_{1})}.\label{eq:ckq}
\end{equation}

Let us now define the space $\left\langle A_{0},A_{1}\right\rangle _{\theta,(r)}$
to consist of all elements $a\in A_{0}+A_{1}$ which can be represented
in the form 
\begin{equation}
a=\sum_{n=-\infty}^{\infty}a_{n}\label{eq:repn}
\end{equation}
for some sequence $\left\{ a_{n}\right\} _{n\in\mathbb{Z}}\in J(\theta,r,A_{0},A_{1})$.
We shall norm this space by 
\begin{equation}
\left\Vert a\right\Vert _{\left\langle A_{0},A_{1}\right\rangle _{\theta,(r)}}:=\inf\left\{ \left\Vert \left\{ a_{n}\right\} _{n\in\mathbb{Z}}\right\Vert _{J(\theta,r,A_{0},A_{1})}\right\} ,\label{eq:few}
\end{equation}
where the infimum is taken over all sequences $\left\{ a_{n}\right\} _{n\in\mathbb{Z}}\in J(\theta,r,A_{0},A_{1})$
which satisfy (\ref{eq:repn}). We can use (\ref{eq:ckq}) to show
that (\ref{eq:few}) defines a norm, rather than merely a seminorm,
and also to show that $\left\langle A_{0},A_{1}\right\rangle _{\theta,(r)}$
is continuously embedded into $A_{0}+A_{1}$. 

If we choose $r=\ln2$, then the above definition of $\left\langle A_{0},A_{1}\right\rangle _{\theta,(r)}$
coincides exactly with the definition of the space denoted by $\left\langle A_{0},A_{1}\right\rangle _{\theta}$
on lines 9--10 of \cite[p.~176]{PeetreJ1971}. If we choose $r=1$
then the same definition of $\left\langle A_{0},A_{1}\right\rangle _{\theta,(r)}$
coincides exactly with the definition of the space, also denoted by
$\left\langle A_{0},A_{1}\right\rangle _{\theta}$, on lines 9--10
of \cite[p.~263]{CwKa}. (Cf.~also \cite[p.~251]{CwikelMKaltonNkmr}.
The paper \cite{CwikelMKaltonNkmr} was catalysed in large part by
ideas appearing in \cite{PeetreJ1971}.) 

It seems intuitively obvious, as was tacitly assumed in \cite{CwKa}
and \cite{CwikelMKaltonNkmr}, that choosing the constant $r$ to
be $\ln2$ or $1$, or indeed any other positive number, gives the
same space $\left\langle A_{0},A_{1}\right\rangle _{\theta}$ to within
equivalence of norms. In the next section, Theorem \ref{thm:EquivPlusMinus},
our main result in this note, will show this explicitly, and will
also provide explicit (but probably not optimal) estimates for the
constants of equivalence. As already mentioned above, this will be
done with the help of an alternative ``continuous'' definition of
these spaces. The existence of this kind of alternative definition
answers one case of a question which was already posed in Probl\`eme
1.2 of \cite[p.~177]{PeetreJ1971}. The much better known real and
complex methods of interpolation have equivalent ``continuous''
and ``discrete'' (or ``periodic'') definitions. (See \cite[pp.~18--19]{LionsJPeetreJ1964}
or \cite[Lemma 3.1.3 p.~41 and Lemma 3.2.3 p.~43]{BerghJLofstromJ1976}
for this result for the real method, and \cite[pp.~1007--1009]{CwikelM1978-complex}
or \cite{AvniE2014} for its analogue for the complex method.) So
we shall see here, not unexpectedly, that the same is true for the
Peetre ``plus-minus'' method%
\footnote{The reason for the name ``plus-minus'', which is sometimes used
for the method for constructing the spaces $\left\langle A_{0},A_{1}\right\rangle _{\theta}$,
comes from a property of (weakly) unconditionally convergent series
which is mentioned, for example, on line 20 of \cite[p.~58]{JansonS1981}.%
}. 

Another potentially useful fact which will emerge from our discussion
here is that the ``continuous'' version of the ``plus-minus''
method, which we shall introduce in Section \ref{sec:CtsDefn} is
in some sense the ``limit'' of its ``discrete'' versions as $r$
tends to $0$. This will be expressed by the formula (\ref{eq:LimCase}).
An analogous result for the complex method is presented in \cite{AvniE2014}.
The analogous result for the real method seems rather obvious. See
Remark \ref{rem:SameForLattices} for another role, perhaps its most
interesting role, which I expect this ``continuous'' definition
of $\left\langle A_{0},A_{1}\right\rangle _{\theta}$ to play.

As I already claimed in a previous version of this note, one can surely
obtain essentially the same result as our main theorem in an alternative
way, by using appropriate modifications of the ideas and methods in
\cite{JansonS1981}, in particular, those on p.~59 of that paper.
I am very grateful to Mieczys\l{}aw Masty\l{}o for subsequently indicating,
in a private communication, how this can indeed be done.

\section{\label{sec:CtsDefn}A ``continuous'' definition of the space $\left\langle A_{0},A_{1}\right\rangle _{\theta}$}

Not surprisingly, in this new definition, sums will be replaced by
integrals. The role previously played by the set of sequences $\Lambda$
will now be played by the set, which we will denote by $\Phi$, of
all functions $\phi:\mathbb{R}\to\mathbb{F}$ which satisfy $\sup_{t\in\mathbb{R}}\left|\phi(t)\right|\le1$
and are locally piecewise continuous, i.e., piecewise continuous on
every bounded interval in $\mathbb{R}$.

For each $\theta\in(0,1)$ the role of the space $J(\theta,r,A_{0},A_{1})$
will now be placed by a space which we will denote by $J(\theta,A_{0},A_{1})$,
or simply by $J$. This will be the space of all functions $u:\mathbb{R}\to A_{0}\cap A_{1}$
which are locally piecewise continuous and \textcolor{blue}{}%
{} for which, for $j\in\left\{ 0,1\right\} $, the improper $A_{j}$
valued Riemann integrals $\int_{-\infty}^{\infty}e^{(j-\theta)t}u(t)dt$
are \textit{unconditionally convergent}, which we will define here
to mean that the integrals $\int_{-R}^{R}e^{(j-\theta)t}\phi(t)u(t)dt$
converge in $A_{j}$ norm as $R\to+\infty$ for every choice of $\phi\in\Phi$.

(For our purposes here we do not need to use Bochner integration or
other more ``advanced'' integration methods for Banach space valued
functions, since we do not need to consider a more general class of
such functions. All integrals appearing here will be ``naive'' Riemann
integrals or improper Riemann integrals of locally piecewise continuous
Banach space valued functions, where the relevant Banach space will
always be $A_{0}\cap A_{1}$. The range of integration for these intervals,
if it is not a bounded or unbounded interval, will be a finite union
of such intervals. The definitions and some of the basic properties
of such integrals can be recalled e.g., by consulting the appendix
of \cite[pp.~257--259]{KatznelsonY1968}.)
\begin{rem}
\label{rem:Cauchy2-1}It will be convenient to explicitly state the
following obvious ``continuous'' analogue of the equivalence of
conditions (i) and (ii) of Remark \ref{rem:CauchyDef}: A function
$u:\mathbb{R}\to A_{0}\cap A_{1}$, which is locally piecewise continuous,
is an element of $J$ if and only if, for $j\in\left\{ 0,1\right\} $
and for each $\varepsilon>0$ and each $\phi\in\Phi$, there exists
a positive number $R(\varepsilon)$ (which could in principle also
depend on $\phi$ and $j$) such that $\left\Vert \int_{R_{1}\le\left|t\right|\le R_{2}}e^{(j-\theta)t}\phi(t)u(t)dt\right\Vert _{A_{j}}\le\varepsilon$
for all numbers $R_{1}$ and $R_{2}$ which satisfy $R(\varepsilon)\le R_{1}<R_{2}$.
\end{rem}

The following lemma is a ``continuous'' analogue of Remark \ref{rem:DiscreteMonotone}.
\begin{lem}
\label{lem:MonotoneCts}Suppose that $U\subset V\subset\mathbb{R}$,
where the each of the sets $U$ and $V$ is the union of finitely
many bounded or unbounded intervals. Suppose that $j\in\left\{ 0,1\right\} $
and that the function $u:\mathbb{R}\to A_{0}\cap A_{1}$ is either

$\mathrm{(i)\,}$ an element of $J$ 

or, alternatively, 

\textup{$\mathrm{(ii)\,}$} is a function which is locally piecewise
continuous and for which 
\begin{equation}
\int_{V}e^{(j-\theta)t}\phi(t)u(t)dt\in A_{j}\mbox{ {\normalcolor \,{\normalcolor for\,\ each}\,}}\phi\in\Phi\,{\normalcolor {\normalcolor \mbox{and\,}}}j\in\left\{ 0,1\right\} .\label{eq:imposed}
\end{equation}

Then, for $j\in\left\{ 0,1\right\} $, 
\begin{equation}
\int_{U}e^{(j-\theta)t}\phi(t)u(t)dt\in A\,\mbox{\,\ for each\,}\phi\in\Phi\label{eq:fef}
\end{equation}
 and 
\begin{equation}
\sup\left\{ \left\Vert \int_{U}e^{(j-\theta)t}\phi(t)u(t)dt\right\Vert _{A{}_{j}}:\,\phi\in\Phi\right\} \le\sup\left\{ \left\Vert \int_{V}e^{(j-\theta)t}\phi(t)u(t)dt\right\Vert _{A{}_{j}}:\,\phi\in\Phi\right\} .\label{eq:MonCts}
\end{equation}

\end{lem}

\noindent \textit{Proof.} For each $\phi\in\Phi$, the form of $U$
and of $V$ ensures that the functions $\phi_{U}:=\phi\chi_{U}$ and
$\phi_{V}:=\phi\chi_{V}$ are also elements of $\Phi$. In case (ii)
the property (\ref{eq:imposed}) is explicitly imposed. In case (i)
the same property follows from the fact that $\int_{V}e^{(j-\theta)t}\phi(t)u(t)dt=\int_{-\infty}^{\infty}e^{(j-\theta)t}\phi_{V}(t)u(t)dt\in A_{j}$
for each $\phi\in\Phi$. This property now implies (\ref{eq:fef})
because $\int_{U}e^{(j-\theta)t}\phi(t)u(t)dt=\int_{V}e^{(j-\theta)t}\phi_{U}(t)u(t)dt$
for each $\phi\in\Phi$. This last formula implies that 
\[
\sup\left\{ \left\Vert \int_{U}e^{(j-\theta)t}\phi(t)u(t)dt\right\Vert _{A{}_{j}}:\,\phi\in\Phi\right\} =\sup\left\{ \left\Vert \int_{V}e^{(j-\theta)t}\phi(t)u(t)dt\right\Vert _{A{}_{j}}:\,\phi\in\Phi_{U}\right\} 
\]
where $\Phi_{U}$ is the set of all functions of the form $\phi_{U}=\phi\chi_{U}$
as $\phi$ ranges over all elements of $\Phi$. Since $\Phi_{U}\subset\Phi$,
we obtain (\ref{eq:MonCts}) and the proof is complete. $\qed$ 

Note that the preceding simple proof has to be valid, and indeed is
valid, even if one or both of the suprema in (\ref{eq:MonCts}) are
infinite, since so far we do not know how to exclude that possibility.
But the next result will imply that both of these suprema are necessarily
finite when $u\in J$.
\begin{fact}
\label{fact:zPrinch}For each $u\in J$, 
\begin{equation}
\lim_{R\to+\infty}\sup\left\{ \left\Vert \int_{\left|t\right|\ge R}e^{(j-\theta)t}\phi(t)u(t)\right\Vert _{A{}_{j}}:\,\phi\in\Phi\right\} =0\mbox{ for \ensuremath{j\in\left\{ 0,1\right\} .}}\label{eq:princh}
\end{equation}

\end{fact}
\noindent \textit{Proof.} (Until we reach the end of this proof,
we still cannot, and do not exclude the possibility that $\sup\left\{ \left\Vert \int_{\left|t\right|\ge R}e^{(j-\theta)t}\phi(t)u(t)dt\right\Vert _{A{}_{j}}:\,\phi\in\Phi\right\} $
could equal $\infty$ for some or all values of $R$.) 

Suppose that (\ref{eq:princh}) does not hold. Then, for at least
one value of $j$ and some $\delta>0$, there exist a strictly increasing
unbounded sequence of positive numbers $\left\{ \rho_{n}\right\} _{n\in\mathbb{N}}$
and a sequence $\left\{ \phi_{n}\right\} _{n\in\mathbb{N}}$ in $\Phi$
such that $\left\Vert \int_{\left|t\right|\ge\rho_{n}}e^{(j-\theta)t}\phi_{n}(t)u(t)dt\right\Vert _{A{}_{j}}\ge\delta$.
For each fixed $n\in\mathbb{N}$, since $\int_{-R}^{R}e^{(j-\theta)t}\phi_{n}(t)u(t)dt$
converges in $A_{j}$ norm as $R\to+\infty$, there exists a number
$R_{n}$ such that $R_{n}>\rho_{n}$ and $\left\Vert \int_{\rho_{n}\le\left|t\right|<R_{n}}e^{(j-\theta)t}\phi_{n}(t)u(t)dt\right\Vert _{A{}_{j}}\ge\delta/2$.
Let us recursively define a sequence $\left\{ n_{k}\right\} _{k\in\mathbb{N}}$
of positive integers such that $n_{1}=1$ and $\rho_{n_{k+1}}>R_{n_{k}}$
for each $k\ge1$. Thus the sets $E_{k}:=\left\{ t\in\mathbb{R}:\rho_{n_{k}}\le\left|t\right|<R_{n_{k}}\right\} $
are pairwise disjoint. This means that the function $\psi=\sum_{k\in\mathbb{N}}\phi_{n_{k}}\chi_{E_{k}}$
is an element of $\Phi$. So (cf.~Remark \ref{rem:Cauchy2-1}) the
quantity $\left\Vert \int_{\rho_{n_{k}}\le\left|t\right|<R_{n_{k}}}e^{(j-\theta)t}\psi(t)u(t)dt\right\Vert _{A_{j}}$
can be made arbitrarily small when $k$ and therefore $\rho_{n_{k}}$
are chosen sufficiently large. But this contradicts the fact that,
for every $k\in\mathbb{N}$, 
\[
\left\Vert \int_{\rho_{n_{k}}\le\left|t\right|<R_{n_{k}}}e^{(j-\theta)t}\psi(t)u(t)dt\right\Vert _{A_{j}}=\left\Vert \int_{\rho_{n_{k}}\le\left|t\right|<R_{n_{k}}}e^{(j-\theta)t}\phi_{n_{k}}(t)u(t)dt\right\Vert _{A_{j}}\ge\frac{\delta}{2}.
\]
 We have therefore proved that (\ref{eq:princh}) holds. $\qed$

Since $\sup\left\{ \left\Vert \int_{-n}^{n}e^{(j-\theta)t}\phi(t)u(t)dt\right\Vert _{A_{j}}:\phi\in\Phi\right\} \le\int_{-n}^{n}e^{(j-\theta)t}\left\Vert u(t)\right\Vert _{A_{j}}dt<\infty$
for each $n\in\mathbb{N}$, we deduce from (\ref{eq:princh}) that
the seminorm 
\begin{equation}
\left\Vert u\right\Vert _{J}:=\sup\left\{ \left\Vert \int_{-\infty}^{\infty}e^{(j-\theta)t}\phi(t)u(t)dt\right\Vert _{A_{j}}:j=0,1,\,\phi\in\Phi\right\} \label{eq:yqw}
\end{equation}
is finite for each $u\in J$. The fact that $\left\Vert \cdot\right\Vert _{J}$
is not quite a norm on $J$ will not cause any difficulties. (One
could of course slightly modify the definition of $J$, e.g., by adding
the requirement that all its elements must be left continuous at all
points, so that $\left\Vert \cdot\right\Vert _{J}$ would become a
norm. But there is no need to do this.) \textcolor{blue}{}%

\begin{rem}
\label{rem:PhiZero}It is clear (analogously to our observation in
Remark \ref{rem:LambdaZero}) that, for each $u\in J$, the value
of the supremum in (\ref{eq:yqw}) will not change if we replace $\Phi$
in (\ref{eq:yqw}) by its subset $\Phi_{0}$ consisting of those of
its functions which have compact support.
\end{rem}

\begin{fact}
\label{fact:Aplus}For each $u\in J$, the limits $\int_{-\infty}^{0}u(t)dt:=\lim_{R\to+\infty}\int_{-R}^{0}u(t)dt$
and $\int_{0}^{\infty}u(t)dt:=\lim_{R\to+\infty}\int_{0}^{R}u(t)dt$
exist, the first with respect to the norm of $A_{0}$, and the second
with respect to the norm of $A_{1}$. Furthermore they satisfy $\left\Vert \int_{-\infty}^{0}u(t)dt\right\Vert _{A_{0}}\le\left\Vert u\right\Vert _{J}$
and $\left\Vert \int_{0}^{\infty}u(t)dt\right\Vert _{A_{1}}\le\left\Vert u\right\Vert _{J}$.
Therefore, for each constant $c\in\mathbb{R}$, the limit $\lim_{R\to+\infty}\int_{-R}^{R+c}u(t)dt$
exists with respect to the norm of $A_{0}+A_{1}$ and does not depend
on $c$. We denote this limit by $\int_{-\infty}^{\infty}u(t)dt$.
It satisfies 
\begin{equation}
\left\Vert \int_{-\infty}^{\infty}u(t)dt\right\Vert _{A_{0}+A_{1}}\le2\left\Vert u\right\Vert _{J}.\label{eq:ejta}
\end{equation}

\end{fact}
\noindent \textit{Proof.} The function $\phi_{0}(t):=\chi_{(-\infty,0]}(t)e^{\theta t}$
is in $\Phi$ and therefore, for each $u\in J$, the integrals $\int_{-R}^{0}u(t)dt=\int_{-R}^{R}e^{-\theta t}\phi_{0}(t)u(t)dt$
converge in $A_{0}$ norm as $R$ tends to $+\infty$ and in fact
\[
\left\Vert \int_{-\infty}^{0}u(t)dt\right\Vert _{A_{0}}=\left\Vert \int_{-\infty}^{\infty}e^{-\theta t}\phi_{0}(t)u(t)dt\right\Vert _{A_{0}}\le\left\Vert u\right\Vert _{J}\,.
\]
Similarly, the function $\phi_{1}(t):=\chi_{[0,\infty)}(t)e^{-(1-\theta)t}$
is in $\Phi$ and therefore, for each $u\in J$, the integrals $\int_{0}^{R}u(t)dt=\int_{-R}^{R}e^{(1-\theta)t}\phi_{1}(t)dt$
converge in $A_{1}$ norm as $R$ tends to $+\infty$ and 
\[
\left\Vert \int_{0}^{\infty}u(t)dt\right\Vert _{A_{1}}=\left\Vert \int_{-\infty}^{\infty}e^{(1-\theta)t}\phi_{1}(t)u(t)dt\right\Vert _{A_{1}}\le\left\Vert u\right\Vert _{J}\,.
\]
The remaining claims in the statement of Fact \ref{fact:Aplus} now
follow obviously and immediately. $\qed$

Let $\left\langle A_{0},A_{1}\right\rangle _{\theta,(0)}$ be the
space of all elements $a\in A_{0}+A_{1}$ which satisfy $a=\int_{-\infty}^{\infty}u(t)dt$
for some $u\in J$. We norm this space by $\left\Vert a\right\Vert _{\left\langle A_{0},A_{1}\right\rangle _{\theta,(0)}}:=\inf\left\Vert u\right\Vert _{J}$
where the infimum is taken over all $u\in J$ for which $a=\int_{-\infty}^{\infty}u(t)dt$.
Analogously to some statements above about the spaces $\left\langle A_{0},A_{1}\right\rangle _{\theta,(r)}$
for $r>0$, here we can use (\ref{eq:ejta}) to show that this infimum
is indeed a norm, rather than merely a seminorm, and also to show
that $\left\langle A_{0},A_{1}\right\rangle _{\theta,(0)}$ is continuously
embedded in $A_{0}+A_{1}$.

We are now ready to state and prove the main theorem of this note.
As an immediate corollary it will provide us with the formula 
\begin{equation}
\left\Vert a\right\Vert _{\left\langle A_{0},A_{1}\right\rangle _{\theta,(0)}}=\lim_{r\searrow0}\left\Vert a\right\Vert _{\left\langle A_{0},A_{1}\right\rangle _{\theta,(r)}}\mbox{ for every }a\in\left\langle A_{0},A_{1}\right\rangle _{\theta,(0)},\label{eq:LimCase}
\end{equation}
which motivates our choice of notation $\left\langle A_{0},A_{1}\right\rangle _{\theta,(0)}$
for the space which we have just introduced.
\begin{thm}
\textup{\label{thm:EquivPlusMinus}For each Banach couple $\left(A_{0},A_{1}\right)$
and each $\theta\in(0,1)$ and each $r>0$, the spaces $\left\langle A_{0},A_{1}\right\rangle _{\theta,(r)}$
and $\left\langle A_{0},A_{1}\right\rangle _{\theta,(0)}$ coincide
to within equivalence of norms. More explicitly, their norms satisfy
\begin{equation}
e^{-r\theta}\left\Vert a\right\Vert _{\left\langle A_{0},A_{1}\right\rangle _{\theta,(r)}}\le\left\Vert a\right\Vert _{\left\langle A_{0},A_{1}\right\rangle _{\theta,(0)}}\le e^{(1-\theta)r}\left\Vert a\right\Vert _{\left\langle A_{0},A_{1}\right\rangle _{\theta,(r)}}\mbox{\,\,\ \ensuremath{\forall}}a\in\left\langle A_{0},A_{1}\right\rangle _{\theta,(0)}.\label{eq:daqw}
\end{equation}
}
\end{thm}
\noindent \textit{Proof.} Our reasoning will be conceptually quite
simple, and quite reminiscent of the straightforward arguments which
can be used (cf.~\cite{LionsJPeetreJ1964} and \cite{BerghJLofstromJ1976})
for establishing the ``classical'' connections between the discrete
and continuous $J$-method constructions of the Lions-Peetre spaces
$\left(A_{0},A_{1}\right)_{\theta,p}$. However here there are rather
more details which have to be carefully checked along the way. 

Suppose first that $a\in\left\langle A_{0},A_{1}\right\rangle _{\theta,(r)}$
for some $r>0$. Given an arbitrary $\varepsilon>0$, choose a sequence
$\left\{ a_{n}\right\} _{n\in\mathbb{Z}}$ in $J(\theta,r,A_{0},A_{1})$,
such that $a=\sum_{n=-\infty}^{\infty}a_{n}$ and 
\[
\left\Vert \left\{ a_{n}\right\} _{n\in\mathbb{Z}}\right\Vert _{J(\theta,r,A_{0},A_{1})}\le\left\Vert a\right\Vert _{\left\langle A_{0},A_{1}\right\rangle _{\theta,(r)}}+\varepsilon.
\]
Now let $u:\mathbb{R}\to A_{0}\cap A_{1}$ be the function defined
by 
\begin{equation}
u=\frac{1}{r}\sum_{n\in\mathbb{Z}}\chi_{[rn,r(n+1))}a_{n}.\label{eq:ItsForm}
\end{equation}
For each $\phi\in\Phi$, for each interval $[rn,r(n+1))$ and for
$j$ equal to either $0$ or $1$ we have that 

\begin{eqnarray}
\int_{rn}^{r(n+1)}e^{(j-\theta)t}\phi(t)u(t)dt & = & \frac{1}{r}e^{(j-\theta)rn}\left(\int_{rn}^{r(n+1)}e^{(j-\theta)(t-rn)}\phi(t)dt\right)a_{n}\nonumber \\
 & = & e^{(j-\theta)rn}\mu_{n}a_{n}\label{eq:rca}
\end{eqnarray}
where $\left|\mu_{n}\right|\le\sup\left\{ \left|e^{(j-\theta)(t-rn)}\phi(t)\right|:t\in[rn,r(n+1)]\right\} $.
So, whether $j=0$ or $j=1$, we obtain that $\left|\mu_{n}\right|\le\max\left\{ e^{(1-\theta)r},1\right\} =e^{(1-\theta)r}$
for all $n$. 

Let $V$ be a set of the form $V=\bigcup_{n\in F}[rn,r(n+1))$ for
some \textit{finite} set $F$ of integers. For $j\in\left\{ 0,1\right\} $
and each $\phi\in\Phi$ we obtain from (\ref{eq:rca}) that 
\begin{equation}
\int_{V}e^{(j-\theta)t}\phi(t)u(t)dt=\sum_{n\in F}e^{(j-\theta)rn}\mu_{n}a_{n}\in A_{j}\label{eq:bwa}
\end{equation}
where, as before, we can assert that the numbers $\mu_{n}$ each satisfy
\begin{equation}
\left|\mu_{n}\right|\le e^{(1-\theta)r}\mbox{ for each }n\in F\mbox{ and for all choices of }\phi\in\Phi.\label{eq:bwaN}
\end{equation}
Let $U$ be an arbitrary union of finitely many arbitrary intervals,
but with the additional condition that it must be a subset of $V$.
The sets $U$ and $V$ satisfy the hypotheses of Lemma \ref{lem:MonotoneCts}.
Furthermore, (\ref{eq:bwa}) shows that $u$ satisfies the condition
(ii) of that lemma. So the lemma justifies the transition from the
first to the second line in the following calculation. The subsequent
transition to the third line follows from (\ref{eq:bwa}) and (\ref{eq:bwaN}).
\begin{eqnarray}
 &  & \sup\left\{ \left\Vert \int_{U}e^{(j-\theta)t}\phi(t)u(t)dt\right\Vert _{A_{j}}:\phi\in\Phi\right\} \nonumber \\
 & \le & \sup\left\{ \left\Vert \int_{V}e^{(j-\theta)t}\phi(t)u(t)dt\right\Vert _{A_{j}}:\phi\in\Phi\right\} \nonumber \\
 & \le & e^{(1-\theta)r}\sup\left\{ \left\Vert \sum_{n\in F}e^{(j-\theta)rn}\lambda_{n}a_{n}\right\Vert _{A_{j}}:\left\{ \lambda_{k}\right\} _{k\in\mathbb{Z}}\in\Lambda\right\} .\label{eq:mywN}
\end{eqnarray}

We shall now use some special cases of these estimates to show that
the function $u$ is an element of $J$:

Let $\rho$ and $R$ be arbitrary numbers which satisfy $r<\rho<R$.
Let $m(\rho)$ be the unique integer for which 
\begin{equation}
rm(\rho)<\rho\le r(m(\rho)+1)\label{eq:gwq}
\end{equation}
 and let $n(R)$ be the unique integer for which $rn(R)>R\ge r(n(R)-1)$.
Then $1\le m(\rho)<n(R)$ and 
\begin{eqnarray*}
U: & = & \left\{ t\in\mathbb{R}:\rho\le\left|t\right|\le R\right\} \subset V:=\left\{ t\in\mathbb{R}:\left|t\right|\in\bigcup_{m(\rho)\le k\le n(R)}[rk,r(k+1)]\right\} \\
 & = & \bigcup_{m(\rho)\le k\le n(R)}\left([-r(k+1),-rk]\cup[rk,r(k+1)]\right)\\
 & = & \bigcup\left\{ [rk,r(k+1)]:k\in F\right\} 
\end{eqnarray*}
where $F$ is the union of the two sets $\left\{ k\in\mathbb{Z}:m(\rho)\le k\le n(R)\right\} $
and 
\[
\left\{ k\in\mathbb{Z}:-n(R)-1\le k\le-m(\rho)-1\right\} .
\]
So, in this case, the inequality provided by (\ref{eq:mywN}) will
justify the first two lines of the following calculation. The transition
to the third line will be justified by Remark \ref{rem:DiscreteMonotone}
together with the fact that $F$ is contained in the set $\left\{ n\in\mathbb{Z}:m(\rho)\le\left|n\right|\le n(R)+1\right\} $.
\begin{eqnarray*}
 &  & \sup\left\{ \left\Vert \int_{\rho\le\left|t\right|\le R}e^{(j-\theta)t}\phi(t)u(t)dt\right\Vert _{A_{j}}:\phi\in\Phi\right\} \\
 & \le & e^{(1-\theta)r}\sup\left\{ \left\Vert \sum_{n\in F}e^{(j-\theta)rn}\lambda_{n}a_{n}\right\Vert _{A_{j}}:\left\{ \lambda_{k}\right\} _{k\in\mathbb{Z}}\in\Lambda\right\} \\
 & \le & e^{(1-\theta)r}\sup\left\{ \left\Vert \sum_{m(\rho)\le\left|n\right|\le n(R)+1}e^{(j-\theta)rn}\lambda_{n}a_{n}\right\Vert _{A_{j}}:\left\{ \lambda_{k}\right\} _{k\in\mathbb{Z}}\in\Lambda\right\} .
\end{eqnarray*}
As an element of $J(\theta,r,A_{0},A_{1})$, the sequence $\left\{ a_{n}\right\} _{n\in\mathbb{Z}}$
must satisfy Condition (iii) of Remark \ref{rem:CauchyDef}. I.e.,
for each $\varepsilon_{0}>0$ there exists an integer $N(\varepsilon_{0})$
such that the expression on the last line of the preceding calculation
is less than $\varepsilon_{0}$ whenever $m(\rho)>N(\varepsilon_{0})$.
Therefore, using (\ref{eq:gwq}), we deduce that the supremum in the
first line of the preceding calculation is less than $\varepsilon_{0}$
whenever $\rho\ge r\left(N(\varepsilon_{0})+1\right)$. Since $\varepsilon_{0}$
is arbitrary this suffices to show (cf.~Remark \ref{rem:Cauchy2-1})
that $u\in J$. (For this we of course also need $u$ to be locally
piecewise continuous. But that is obvious from (\ref{eq:ItsForm}).)

The fact that $u\in J$ implies in turn (cf.~Fact \ref{fact:Aplus})
that the integral $\int_{-\infty}^{\infty}u(t)dt$ exists as an element
of $A_{0}+A_{1}$ equal to the limit in $A_{0}+A_{1}$ norm as $R\to+\infty$
of the integrals $\int_{-R}^{R}u(t)dt$. This limit must coincide
with 
\[
\lim_{n\to\infty}\int_{-rn}^{rn}u(t)dt=\lim_{n\to\infty}\sum_{k=-n}^{n-1}\int_{rk}^{r(k+1)}u(t)dt=\lim_{n\to\infty}\sum_{k=-n}^{n-1}a_{k}=\lim_{n\to\infty}\sum_{k=-n}^{n}a_{k}=a.
\]
(Here we have used the result described just before (\ref{eq:ckq})
and chosen the constant $k_{0}$ appearing there to equal $-1$.)
We also claim that 
\begin{equation}
\left\Vert u\right\Vert _{J}\le e^{(1-\theta)r}\left\Vert \left\{ a_{n}\right\} _{n\in\mathbb{Z}}\right\Vert _{J\left(\theta,r,\left(A_{0},A_{1}\right)\right)}.\label{eq:pgr}
\end{equation}
To show this, it will suffice (cf.~Remark \ref{rem:PhiZero}) to
show that 
\begin{equation}
\sup\left\{ \left\Vert \int_{U}e^{(j-\theta)t}\phi(t)u(t)dt\right\Vert _{A_{j}}:\phi\in\Phi\right\} \le e^{(1-\theta)r}\left\Vert \left\{ a_{n}\right\} _{n\in\mathbb{Z}}\right\Vert _{J\left(\theta,r,\left(A_{0},A_{1}\right)\right)}\label{eq:bfq}
\end{equation}
for $j\in\left\{ 0,1\right\} $ and for each bounded interval $U\subset\mathbb{R}$.
Each such $U$ can of course be contained in a set of the form $V=\bigcup_{n\in F}[rn,r(n+1))$
for some finite set $F$ of integers. So we can apply the inequalities
(\ref{eq:mywN}) together with the fact (cf.~Remark \ref{rem:DiscreteMonotone})
that
\begin{eqnarray*}
 &  & \sup\left\{ \left\Vert \sum_{n\in F}e^{(j-\theta)rn}\lambda_{n}a_{n}\right\Vert _{A_{j}}:\left\{ \lambda_{k}\right\} _{k\in\mathbb{Z}}\in\Lambda\right\} \\
 & \le & \sup\left\{ \left\Vert \sum_{n\in\mathbb{Z}}e^{(j-\theta)rn}\lambda_{n}a_{n}\right\Vert _{A_{j}}:\left\{ \lambda_{k}\right\} _{k\in\mathbb{Z}}\in\Lambda\right\} \le\left\Vert \left\{ a_{n}\right\} _{n\in\mathbb{Z}}\right\Vert _{J\left(\theta,r,\left(A_{0},A_{1}\right)\right)}
\end{eqnarray*}
to obtain (\ref{eq:bfq}) for each such $U$ and therefore also (\ref{eq:pgr}).
Thus we have shown that $a\in\left\langle A_{0},A_{1}\right\rangle _{\theta,(0)}$
and $\left\Vert a\right\Vert _{\left\langle A_{0},A_{1}\right\rangle _{\theta,(0)}}\le\left\Vert u\right\Vert _{J}\le e^{(1-\theta)r}\left(\left\Vert a\right\Vert _{\left\langle A_{0},A_{1}\right\rangle _{\theta,(r)}}+\varepsilon\right)$.
Since $a$ is an arbitrary element of $\left\langle A_{0},A_{1}\right\rangle _{\theta,(r)}$
and since $\varepsilon$ is an arbitrary positive number, we conclude
that $\left\langle A_{0},A_{1}\right\rangle _{\theta,(r)}$ is continuously
embedded in $\left\langle A_{0},A_{1}\right\rangle _{\theta,(0)}$
and that the norms of these two spaces satisfy the second inequality
in (\ref{eq:daqw}).

Now suppose, conversely, that $a\in\left\langle A_{0},A_{1}\right\rangle _{\theta,(0)}$.
Given an arbitrary $\varepsilon>0$, let $u$ be an element of $J$
for which $a=\int_{-\infty}^{\infty}u(t)dt$ and $\left\Vert u\right\Vert _{J}\le\left\Vert a\right\Vert _{\left\langle A_{0},A_{1}\right\rangle _{\theta,(0)}}+\varepsilon$.
Define the sequence $\left\{ a_{n}\right\} _{n\in\mathbb{Z}}$ of
elements of $A_{0}\cap A_{1}$ by $a_{n}=\int_{rn}^{r(n+1)}u(t)dt$
for each $n\in\mathbb{Z}$. Given an arbitrary sequence $\left\{ \lambda_{k}\right\} _{k\in\mathbb{Z}}$
in $\Lambda$, define a function $\psi\in\Phi$ associated with $\left\{ \lambda_{k}\right\} _{k\in\mathbb{Z}}$
by setting $\psi=\sum_{k\in\mathbb{Z}}\lambda_{k}\chi_{[rk,r(k+1))}$.
Then, for every finite subset $F$ of $\mathbb{Z}$, and for $j\in\left\{ 0,1\right\} $.
we have 
\begin{eqnarray*}
\sum_{n\in F}e^{(j-\theta)rn}\lambda_{n}a_{n} & = & \sum_{n\in F}\int_{rn}^{r(n+1)}e^{(j-\theta)rn}\psi(t)u(t)dt\\
 & = & \sum_{n\in F}\int_{rn}^{r(n+1)}e^{(j-\theta)t}\xi_{j}(t)u(t)dt\,,
\end{eqnarray*}
where $\xi_{j}(t)=\sum_{n\in F}e^{(j-\theta)\left(rn-t\right)}\psi(t)$$\chi_{[rn,r(n+1))}$.
It is clear that the function $\xi_{j}$ is locally piecewise continuous
and that, for all $t\in\mathbb{R}$, we have $\left|\xi_{0}(t)\right|\le e^{r\theta}$
and $\left|\xi_{1}(t)\right|\le1$. So $e^{-r\theta}\xi_{j}\in\Phi$
for both values of $j$. Therefore, for every $\left\{ \lambda_{k}\right\} _{k\in\mathbb{Z}}$
in $\Lambda$ and every finite set $F\subset\mathbb{Z},$ we have
\begin{equation}
\left\Vert \sum_{n\in F}e^{(j-\theta)rn}\lambda_{n}a_{n}\right\Vert _{A_{j}}\le e^{r\theta}\sup\left\{ \left\Vert \int_{E_{F}}e^{(j-\theta)t}\phi(t)u(t)dt\right\Vert _{A_{j}}:\phi\in\Phi\right\} ,\label{eq:mbdma}
\end{equation}
where 
\begin{equation}
E_{F}=\bigcup_{n\in F}[rn,r(n+1)).\label{eq:deef}
\end{equation}
Let us now specialize to the case where $F$ is a set of the form
$F=\left\{ n\in\mathbb{Z}:n_{1}\le\left|n\right|\le n_{2}\right\} $
for integers $n_{1}$ and $n_{2}$ which satisfy $n_{2}>n_{1}>2$.
Then $E_{F}$ is contained in $\left\{ t\in\mathbb{R}:\left|t\right|\ge R_{1}\right\} $
whenever the number $R_{1}$ satisfies 
\begin{equation}
0<R_{1}<r(n_{1}-1),\label{eq:rnetc}
\end{equation}
and, since $u\in J$, we can apply Lemma \ref{lem:MonotoneCts} (case
(i)) with $U=E_{F}$ and $V=\left\{ t\in\mathbb{R}:R_{1}\le\left|t\right|\right\} $
to deduce from (\ref{eq:mbdma}) that 
\begin{equation}
\left\Vert \sum_{n_{1}\le\left|n\right|\le n_{2}}e^{(j-\theta)rn}\lambda_{n}a_{n}\right\Vert _{A_{j}}\le e^{r\theta}\sup\left\{ \left\Vert \int_{\left|t\right|\ge R_{1}}e^{(j-\theta)t}\phi(t)u(t)dt\right\Vert _{A_{j}}:\phi\in\Phi\right\} \label{eq:bqe}
\end{equation}

In view of Fact \ref{fact:zPrinch}, the right side of (\ref{eq:bqe})
can be made arbitrarily small by choosing $R_{1}$ sufficiently large.
So, in view of (\ref{eq:rnetc}), the left side of (\ref{eq:bqe})
is arbitrarily small whenever $n_{1}$ is sufficiently large. Since
this property holds for $j=0,1$ and for each $\left\{ \lambda_{k}\right\} _{k\in\mathbb{Z}}$
in $\Lambda$, we see (cf.~Remark \ref{rem:CauchyDef}) that the
sequence $\left\{ a_{n}\right\} _{n\in\mathbb{Z}}$ is an element
of $J(\theta,r,A_{0},A_{1})$. Therefore the sum $\sum_{n=-\infty}^{\infty}a_{n}$
is an element of $\left\langle A_{0},A_{1}\right\rangle _{\theta,(r)}$
which is the limit in $A_{0}+A_{1}$ as $N$ tends to $\infty$ of
the partial sums 
\[
\sum_{n=-N}^{N}a_{n}=\sum_{n=-N}^{N}\int_{rn}^{r(n+1)}u(t)dt=\int_{-rN}^{r(N+1)}u(t)dt.
\]
So this limit equals $\int_{-\infty}^{\infty}u(t)dt=a$. (Here we
have used Fact \ref{fact:Aplus}, with the constant $c$ appearing
there now chosen to equal $r$.) 

In order to estimate the norm of $a$ in $\left\langle A_{0},A_{1}\right\rangle _{\theta,(r)}$,
we once again apply Lemma \ref{lem:MonotoneCts} (case (i)) to (\ref{eq:mbdma})
and (\ref{eq:deef}), but this time with $U=E_{F}$ for an \textit{arbitrary}
finite subset $F$ of $\mathbb{Z}$ and with $V=\mathbb{R}$. This
gives that $\left\Vert \sum_{n\in F}e^{(j-\theta)rn}\lambda_{n}a_{n}\right\Vert _{A_{j}}\le e^{r\theta}\left\Vert u\right\Vert _{J}$
for $j\in\left\{ 0,1\right\} $ and every $\left\{ \lambda_{k}\right\} _{k\in\mathbb{Z}}\in\Lambda$,
independently of our choice of $F$. Therefore (cf.~Remark \ref{rem:LambdaZero})
we have that $\left\Vert \left\{ a_{n}\right\} _{n\in\mathbb{Z}}\right\Vert _{J(\theta,r,A_{0},A_{1})}\le e^{r\theta}\left\Vert u\right\Vert _{J}$
and, consequently, 
\[
\left\Vert a\right\Vert _{\left\langle A_{0},A_{1}\right\rangle _{\theta,(r)}}\le e^{r\theta}\left(\left\Vert a\right\Vert _{\left\langle A_{0},A_{1}\right\rangle _{\theta,(0)}}+\varepsilon\right).
\]
This completes the proof that $\left\langle A_{0},A_{1}\right\rangle _{\theta,(0)}$
is continuously embedded in $\left\langle A_{0},A_{1}\right\rangle _{\theta,(r)}$
and that the norms of these two spaces satisfy the first inequality
in (\ref{eq:daqw}). Therefore this also completes the proof of the
theorem. $\qed$ 
\begin{rem}
Of course the inequalities (\ref{eq:daqw}) can be used to obtain
estimates for the equivalence constants between the norms $\left\Vert \cdot\right\Vert _{\left\langle A_{0},A_{1}\right\rangle _{\theta,(r_{1})}}$
and $\left\Vert \cdot\right\Vert _{\left\langle A_{0},A_{1}\right\rangle _{\theta,(r_{2})}}$
for any two positive numbers $r_{1}$ and $r_{2}$, should anyone
ever need them. But these are very unlikely to be optimal estimates,
given that they are obtained via passage through the norm $\left\Vert \cdot\right\Vert _{\left\langle A_{0},A_{1}\right\rangle _{\theta,(0)}}$.
One can probably obtain better estimates fairly easily in the case
where $r_{1}=2r_{2}$ or $r_{1}$ is some other integer multiple of
$r_{2}$. 
\end{rem}

\section{\label{sec:PeetreInCalderon}The inclusion $\left\langle A_{0},A_{1}\right\rangle _{\theta}\subset\left[A_{0},A_{1}\right]_{\theta}$}

Already in 1971 (see the formula (1.2) on p.~176 of \cite{PeetreJ1971}),
Jaak Peetre remarked that the space $\left\langle A_{0},A_{1}\right\rangle _{\theta}$
is continuously embedded in the complex interpolation space $\left[A_{0},A_{1}\right]_{\theta}$
introduced by Alberto Calder\'on in \cite{CalderonA1964}. (In fact
the remarks in \cite{PeetreJ1971} show that $\left\langle A_{0},A_{1}\right\rangle _{\theta}$
is contained in a possibly smaller ``periodic/discrete'' variant
of $\left[A_{0},A_{1}\right]_{\theta}$ which would later be shown
(see \cite{CwikelM1978-complex}) to coincide with $\left[A_{0},A_{1}\right]_{\theta}$.)
For the reader's convenience, in Theorem \ref{thm:PeetreInCalderon}
of this section, we explicitly present the simple proof that $\left\langle A_{0},A_{1}\right\rangle _{\theta}\subset\left[A_{0},A_{1}\right]_{\theta}$.
We make a point of noting that this inclusion is a norm one embedding,
and we also briefly mention another norm one estimate in terms of
``periodic'' norms on $\left[A_{0},A_{1}\right]_{\theta}$. 
\begin{rem}
\label{rem:SameForLattices}The inclusion $\left\langle A_{0},A_{1}\right\rangle _{\theta}\subset\left[A_{0},A_{1}\right]_{\theta}$
is the easy part of the proof of another interesting property of the
space $\left\langle A_{0},A_{1}\right\rangle _{\theta}$, namely that
it \textit{coincides} with $\left[A_{0},A_{1}\right]_{\theta}$, to
within equivalence of norms, whenever $\left(A_{0},A_{1}\right)$
is a couple of complexified Banach lattices of measurable functions
on the same underlying measure space. This property has sometimes
proved useful, for example in \cite{CwKa}, and I hope to use it further
in the very near future, in a forthcoming paper (which in fact has
been the main motivation for me to write this note). I also conjecture
that this coincidence to within equivalence of norms of $\left\langle A_{0},A_{1}\right\rangle _{\theta}$
and $\left[A_{0},A_{1}\right]_{\theta}$ for couples of lattices is
even an \textit{isometry} when $\left\langle A_{0},A_{1}\right\rangle _{\theta}$
is equipped with the norm of $\left\langle A_{0},A_{1}\right\rangle _{\theta,(0)}$,
i.e., the norm associated with the apparently new ``continuous''
method for its construction presented in Section \ref{sec:CtsDefn}.
Hence my interest in explicitly noting that the embedding $\left\langle A_{0},A_{1}\right\rangle _{\theta,(r)}\subset\left[A_{0},A_{1}\right]_{\theta}$
has norm one.
\end{rem}

\begin{rem}
\label{rem:SvanteExample}The contents of the preceding remark make
it almost compulsory to mention an example due to Svante Janson \cite[Example 6, p.~62]{JansonS1981}.
This example can be used to show that $\left\langle A_{0},A_{1}\right\rangle _{\theta}\ne\left[A_{0},A_{1}\right]_{\theta}$
for some Banach couples $\left(A_{0},A_{1}\right)$ which do not have
the above mentioned lattice structure. More explicitly, if $\left(A_{0},A_{1}\right)$
is the couple $\left(FL_{0}^{1},FL_{1}^{1}\right)$ (in the notation
of \cite[p.~81]{CwikelMJansonS2007}), then the above-mentioned Example
6, combined with Theorems 3 and 5 on pages 57 and 59 of \cite{JansonS1981}
shows that the sequence space $\left\langle FL_{0}^{1},FL_{1}^{1}\right\rangle _{\theta}$
is contained in the weighted $\ell^{2}$ space $\left\{ \left\{ \lambda_{n}\right\} _{n\in\mathbb{Z}}:\sum_{n\in\mathbb{Z}}\left|e^{\theta n}\lambda_{n}\right|^{2}<\infty\right\} $.
On the other hand, as shown in the proof of Theorem 22 on p.~68 of
\cite{JansonS1981}, $\left[FL_{0}^{1},FL_{1}^{1}\right]_{\theta}$
contains all complex sequences $\left\{ \lambda_{n}\right\} _{n\in\mathbb{Z}}$
for which $\left\{ e^{\theta n}\lambda_{n}\right\} _{n\in\mathbb{Z}}$
is the sequence of Fourier coefficients of some function in $L^{1}(\mathbb{T})$.
Since $L^{2}(\mathbb{T})$ is strictly smaller than $L^{1}(\mathbb{T})$
this shows that $\left\langle FL_{0}^{1},FL_{1}^{1}\right\rangle _{\theta}$
is strictly smaller than $\left[FL_{0}^{1},FL_{1}^{1}\right]_{\theta}$.\end{rem}
\begin{thm}
\label{thm:PeetreInCalderon}$($Cf.~\cite[(1.2) p.~176]{PeetreJ1971}
and \cite{JansonS1981}.$)$ For every Banach couple $\left(A_{0},A_{1}\right)$
of complex Banach spaces, and for each $r\ge0$ and $\theta\in(0,1)$,
the inclusion 
\[
\left\langle A_{0},A_{1}\right\rangle _{\theta,(r)}\subset\left[A_{0},A_{1}\right]_{\theta}
\]
holds, and 
\begin{equation}
\left\Vert b\right\Vert _{\left[A_{0},A_{1}\right]_{\theta}}\le\left\Vert b\right\Vert _{\left\langle A_{0},A_{1}\right\rangle _{\theta.(r)}}\mbox{ for every\,\,}b\in\left\langle A_{0},A_{1}\right\rangle _{\theta,(r)}.\label{eq:GoodEstimate}
\end{equation}

\end{thm}
\noindent \textit{Proof.} It suffices to treat the case where $r>0$,
since the case where $r=0$ will then follow immediately from Theorem
\ref{thm:EquivPlusMinus} and (\ref{eq:LimCase}). 

Given an arbitrary sequence $\left\{ a_{n}\right\} _{n\in\mathbb{Z}}$
of elements of $A_{0}\cap A_{1}$, let $b_{U}=\sum_{n\in U}a_{n}$
for every finite subset $U$ of $\mathbb{Z}$. Obviously $b_{U}$
is an element of $\left[A_{0},A_{1}\right]_{\theta}$ since it is
an element of $A_{0}\cap A_{1}$. But we want to control its norm
in $\left[A_{0},A_{1}\right]_{\theta}$. For each $\delta>0$, we
introduce the function $f_{\delta}(z)=e^{\delta(z-\theta)^{2}}\sum_{n\in U}e^{(z-\theta)rn}a_{n}$
which is clearly an element of Calder\'on's space $\mathcal{F}\left(A_{0},A_{1}\right)$
and satisfies $f_{\delta}(\theta)=b_{U}$ and also 
\begin{eqnarray*}
\left\Vert f_{\delta}\right\Vert _{\mathcal{F}\left(A_{0},A_{1}\right)} & = & \sup\left\{ \left\Vert f_{\delta}(j+it)\right\Vert _{A_{j}}:j\in\left\{ 0,1\right\} ,\, t\in\mathbb{R}\right\} \\
 & \le & e^{\delta}\sup\left\{ \left\Vert \sum_{n\in U}e^{(j+it-\theta)rn}a_{n}\right\Vert _{A_{j}}:j\in\left\{ 0,1\right\} ,\, t\in\mathbb{R}\right\} \\
 & \le & e^{\delta}\sup\left\{ \left\Vert \sum_{n\in U}\lambda_{n}e^{(j-\theta)rn}a_{n}\right\Vert _{A_{j}}:j\in\left\{ 0,1\right\} ,\,\left\{ \lambda_{k}\right\} _{k\in\mathbb{Z}}\in\Lambda\right\} .
\end{eqnarray*}
Since $\delta$ can be chosen arbitrarily small, we deduce that 
\begin{equation}
\left\Vert b_{U}\right\Vert _{\left[A_{0},A_{1}\right]_{\theta}}\le\sup\left\{ \left\Vert \sum_{n\in U}\lambda_{n}e^{(j-\theta)rn}a_{n}\right\Vert _{A_{j}}:j\in\left\{ 0,1\right\} ,\,\left\{ \lambda_{k}\right\} _{k\in\mathbb{Z}}\in\Lambda\right\} .\label{eq:blink}
\end{equation}

Let $b$ be an arbitrary element of $\left\langle A_{0},A_{1}\right\rangle _{\theta,(r)}$
and let $\varepsilon$ be an arbitrary positive number. Then there
exists a sequence $\left\{ a_{n}\right\} _{n\in\mathbb{Z}}$ in $J\left(\theta,r,A_{0},A_{1}\right)$
with $\left\Vert \left\{ a_{n}\right\} _{n\in\mathbb{Z}}\right\Vert _{J(\theta,r,A_{0},A_{1})}\le\left(1+\varepsilon\right)\left\Vert b\right\Vert _{\left\langle A_{0},A_{1}\right\rangle _{\theta,(r)}}$
and such that $b=\sum_{n\in\mathbb{Z}}a_{n}$ (convergence in $A_{0}+A_{1}$).
As above, let us define $b_{U}=\sum_{n\in U}a_{n}$ for every finite
subset $U$ of $\mathbb{Z}$. Since $\left\{ a_{n}\right\} _{n\in\mathbb{Z}}$
necessarily satisfies Condition (iii) of Remark \ref{rem:CauchyDef},
we can deduce from (\ref{eq:blink}) that, for each $\varepsilon>0$,
there exists a positive integer $N(\varepsilon)$ such that $\left\Vert b_{U}\right\Vert _{\left[A_{0},A_{1}\right]_{\theta}}<\varepsilon$
whenever $U\subset\left\{ n\in\mathbb{Z}:\left|n\right|\ge N(\varepsilon)\right\} $.
This shows that the sequence $\left\{ s_{k}\right\} _{k\in\mathbb{N}}$
defined by $s_{k}=\sum_{n=-k}^{k}a_{n}$ is a Cauchy sequence in $\left[A_{0},A_{1}\right]_{\theta}$.
Since this sequence converges in $A_{0}+A_{1}$ to $b$, this means
that $b\in\left[A_{0},A_{1}\right]_{\theta}$ and also that $\left\Vert b\right\Vert _{\left[A_{0},A_{1}\right]_{\theta}}=\lim_{k\to\infty}\left\Vert s_{k}\right\Vert _{\left[A_{0},A_{1}\right]_{\theta}}$.
Letting $U_{k}$ be the set of integers $n$ satisfying $\left|n\right|\le k$,
we use (\ref{eq:blink}) again, and then (\ref{eq:itijgts}), with
$U=U_{k}$ and $V=\mathbb{Z}$, to obtain that
\[
\left\Vert s_{k}\right\Vert _{\left[A_{0},A_{1}\right]_{\theta}}=\left\Vert b_{U_{k}}\right\Vert _{\left[A_{0},A_{1}\right]_{\theta}}\le\left\Vert \left\{ a_{n}\right\} _{n\in\mathbb{Z}}\right\Vert _{J(\theta,r,A_{0},A_{1})}\le\left(1+\varepsilon\right)\left\Vert b\right\Vert _{\left\langle A_{0},A_{1}\right\rangle _{\theta,(r)}}.
\]
Since $\varepsilon$ is arbitrary, this completes the proof. $\qed$ 
\begin{rem}
For each finite set $U\subset\mathbb{Z}$, the function $f(z)=\sum_{n\in U}e^{(z-\theta)rn}a_{n}$
satisfies $f(z+2\pi i/r)$ for all $z\in\mathbb{C}$. Therefore, if
we choose $\lambda=2\pi/r$, this function is an element of the ``periodic''
space $\mathcal{F}_{\lambda}\left(A_{0},A_{1}\right)$ as defined
on p.~1007 of \cite{CwikelM1978-complex}. (Cf.~also \cite[pp.~161--162]{AvniE2014}.)
Therefore, for each $r>0$, the same reasoning as in the proof of
the previous theorem in fact gives a slightly stronger estimate than
(\ref{eq:GoodEstimate}), since it is expressed in terms of the larger
``periodic'' norm $\left\Vert \cdot\right\Vert _{\left[A_{0},A_{1}\right]_{\theta}^{\lambda}}$,
defined on that same page of \cite{CwikelM1978-complex}. I.e., we
obtain that 
\[
\left\Vert b\right\Vert _{\left[A_{0},A_{1}\right]_{\theta}^{2\pi/r}}\le\left\Vert b\right\Vert _{\left\langle A_{0},A_{1}\right\rangle _{\theta,(r)}}\mbox{ for all\,}r>0,\,\theta\in(0,1)\mbox{ and\,}b\in\left\langle A_{0},A_{1}\right\rangle _{\theta,(r)}.
\]
 
\end{rem}

\section{\label{sec:Appendix}Appendix - Completeness of the spaces $J(\theta,r,A_{0},A_{1})$
and $\left\langle A_{0},A_{1}\right\rangle _{\theta}$}

For the convenience of those readers who may happen to be less familiar
with these kinds of topics, here is a proof, via more or less routine
arguments, of the completeness of $J(\theta,r,A_{0},A_{1})$ and,
consequently, of $\left\langle A_{0},A_{1}\right\rangle _{\theta(r)}$
for each $r>0$ and $\theta\in(0,1)$. The completeness of $\left\langle A_{0},A_{1}\right\rangle _{\theta,(r)}$
for $r=0$ can then be deduced immediately from Theorem \ref{thm:EquivPlusMinus}.
(So we see that the obvious fact that the space $J$ is not complete
nor even normed, does not create any difficulties.) Naturally enough,
the fact that $\left\langle A_{0},A_{1}\right\rangle _{\theta}$ is
complete was noticed and mentioned already in \cite{PeetreJ1971}. 

We begin by considering two-sided sequences $\left\{ a_{n}\right\} _{n\in\mathbb{Z}}$
which take values in an arbitrary Banach space $A$. As in Example
2.4 of \cite[p.~247]{CwikelMKaltonNkmr}, we let $UC(A)$ denote the
space of all such sequences $\left\{ a_{n}\right\} _{n\in\mathbb{Z}}$
for which $\sum_{n\in\mathbb{Z}}a_{n}$ is unconditionally convergent
in $A$. As already remarked on p.~\pageref{pCzero}, an $A$ valued
sequence $\left\{ a_{n}\right\} _{n\in\mathbb{Z}}$ is in $UC(A)$
if and only if it satisfies the condition (\ref{eq:LikeCzero}). Since
$A$ is complete, (\ref{eq:LikeCzero}) is equivalent, in turn, to
\begin{equation}
\lim_{N\to\infty}\sup\left\{ \left\Vert \sum_{\left|n\right|\ge N}\lambda_{n}a_{n}\right\Vert _{A}:\left\{ \lambda_{k}\right\} _{k\in\mathbb{Z}}\in\Lambda_{0}\right\} =0.\label{eq:bwq}
\end{equation}
Condition (\ref{eq:bwq}) ensures that the functional 
\[
\left\Vert \left\{ a_{n}\right\} _{n\in\mathbb{Z}}\right\Vert _{UC(A)}:=\sup\left\{ \left\Vert \sum_{n\in\mathbb{Z}}\lambda_{n}a_{n}\right\Vert _{A}:\left\{ \lambda_{k}\right\} _{k\in\mathbb{Z}}\in\Lambda_{0}\right\} 
\]
is finite for each $\left\{ a_{n}\right\} _{n\in\mathbb{Z}}\in UC(A)$
and therefore defines a norm%
\footnote{This norm is different from, but equivalent to the norm chosen for
this space in \cite{CwikelMKaltonNkmr}. %
} on $UC(A)$. Note that this norm satisfies 
\begin{equation}
\left\Vert a_{n_{0}}\right\Vert _{A}\le\left\Vert \left\{ a_{n}\right\} _{n\in\mathbb{Z}}\right\Vert _{UC(A)}\mbox{\,\ for each fixed\,}n_{0}\in\mathbb{Z}\mbox{ and each\,}\left\{ a_{n}\right\} _{n\in\mathbb{Z}}\in UC(A).\label{eq:ptw}
\end{equation}

The main step of our proof will be to show that $UC$$(A)$ is complete
with respect to this norm. (This will confirm, as asserted on p.~247
of \cite{CwikelMKaltonNkmr}, that the mapping $UC$ is indeed a ``pseudolattice''
(cf.~Definition 2.1 of \cite[p.~246]{CwikelMKaltonNkmr}).)

For each sequence $a=\left\{ a_{n}\right\} _{n\in\mathbb{Z}}$ in
$UC(A)$ let $Ta$ be the numerical sequence $\left\{ \left(Ta\right)_{N}\right\} _{N\in\mathbb{Z}}$
defined by 
\[
(Ta)_{N}=\sup\left\{ \left\Vert \sum_{\left|n\right|\ge N}\lambda_{n}a_{n}\right\Vert _{A}:\left\{ \lambda_{k}\right\} _{k\in\mathbb{Z}}\in\Lambda_{0}\right\} .
\]
This defines a (nonlinear) map $T:UC(A)\to c_{0}$ which satisfies
\begin{equation}
\left\Vert Ta-Tb\right\Vert _{\ell^{\infty}}\le\left\Vert a-b\right\Vert _{UC(A)}\label{eq:Lipschitz}
\end{equation}
for all pairs of elements $a=\left\{ a_{n}\right\} _{n\in\mathbb{Z}}$
and $b=\left\{ b_{n}\right\} _{n\in\mathbb{Z}}$ in $UC(A)$. Now
suppose that $\left\{ a_{m}\right\} _{m\in\mathbb{N}}$ is a Cauchy
sequence in $UC(A)$. For each fixed $m\in\mathbb{N}$, the element
$a_{m}$ is a two-sided $A$ valued sequence which we will denote
by $\left\{ a_{m,n}\right\} _{n\in\mathbb{Z}}$ and $Ta_{m}$ is a
one-sided numerical sequence which we shall denote by $\left\{ (Ta_{m})_{k}\right\} _{k\in\mathbb{N}}$.
It follows from (\ref{eq:Lipschitz}) that the sequence $\left\{ Ta_{m}\right\} _{m\in\mathbb{N}}$
is a Cauchy sequence in $c_{0}$. Therefore the same kinds of standard
arguments used to characterize the relatively compact subsets of $c_{0}$
show that the numbers $\rho_{N}:=\sup\left\{ \left|\left(Ta_{m}\right)_{k}\right|:k\ge N,\, m\in\mathbb{N}\right\} $
satisfy $\lim_{N\to\infty}\rho_{N}=0$. 

In view of (\ref{eq:ptw}), for each fixed $n_{0}\in\mathbb{Z}$,
the sequence $\left\{ a_{m,n_{0}}\right\} _{m\in\mathbb{N}}$ is Cauchy
in $A$ and therefore there exists an $A$ valued sequence $b=\left\{ b_{n}\right\} _{n\in\mathbb{Z}}$
such that $\lim_{m\to\infty}\left\Vert a_{m,n}-b_{n}\right\Vert _{A}=0$
for each fixed $n\in\mathbb{Z}$. For each $N\in\mathbb{N}$ and each
$\left\{ \lambda_{k}\right\} _{k\in\mathbb{Z}}\in\Lambda_{0}$, we
have $\left\Vert \sum_{\left|n\right|\ge N}\lambda_{n}b_{n}\right\Vert _{A}=\lim_{m\to\infty}\left\Vert \sum_{\left|n\right|\ge N}\lambda_{n}a_{m,n}\right\Vert _{A}\le\rho_{N}$.
This shows see that the sequence $b$ has the property (\ref{eq:bwq})
and is therefore an element of $UC(A)$. 

Given any $\varepsilon>0$, let us choose $N=N(\varepsilon)$ sufficiently
large so that $\rho_{N}<\epsilon/4$. Then, choose $M=M(\varepsilon)$
such that $\left\Vert a_{m,n}-b_{n}\right\Vert _{A}\le\epsilon/4N$
for each integer $m\ge M$ and each integer $n\in\left[-N,N\right]$.
These choices of $N$ and $M$ ensure that, for each $m\ge M$, and
each $\left\{ \lambda_{k}\right\} _{k\in\mathbb{Z}}$ in $\Lambda_{0}$,
\begin{eqnarray*}
\left\Vert \sum_{n\in\mathbb{Z}}\lambda_{n}(a_{m,n}-b_{n})\right\Vert _{A} & \le & \left\Vert \sum_{n=-N+1}^{N-1}\lambda_{n}(a_{m,n}-b_{n})\right\Vert _{A}+\left\Vert \sum_{\left|n\right|\ge N}\lambda_{n}(a_{m,n}-b_{n})\right\Vert _{A}\\
 & \le & \sum_{n=-N+1}^{N-1}\left\Vert a_{m,n}-b_{n}\right\Vert +\left\Vert \sum_{\left|n\right|\ge N}\lambda_{n}a_{m,n}\right\Vert _{A}+\left\Vert \sum_{\left|n\right|\ge N}\lambda_{n}b_{n}\right\Vert _{A}\\
 & \le & (2N-1)\frac{\varepsilon}{4N}+\rho_{N}+\rho_{N}<\varepsilon.
\end{eqnarray*}
Taking the supremum as $\left\{ \lambda_{k}\right\} _{k\in\mathbb{Z}}$
ranges over all of $\Lambda_{0}$, we obtain that $\left\Vert b-a_{m}\right\Vert _{UC(A)}\le\varepsilon$
for all $m\ge M$. This completes the proof that $UC(A)$ is complete.

The rest of the proof of the completeness of $J(\theta,r,A_{0},A_{1})$
and then of $\left\langle A_{0},A_{1}\right\rangle _{\theta}$ is
a special case of the argument indicated immediately after Definition
2.11 on p.~249 of \cite{CwikelMKaltonNkmr}. More explicitly, we
now observe that the space $J(\theta,r,A_{0},A_{1})$ is the intersection
of a ``weighted'' version of $UC(A_{0})$ with a differently ``weighted''
version of $UC(A_{1})$. So the completeness of $J(\theta,r,A_{0},A_{1})$
is a simple consequence of the completeness of each of the two spaces
$UC\left(A_{0}\right)$ and $UC\left(A_{1}\right)$ and the fact that
$A_{0}$ and $A_{1}$ are both continuously embedded in the Banach
space $A_{0}+A_{1}$. The completeness of $\left\langle A_{0},A_{1}\right\rangle _{\theta}$
can be very easily and ``routinely'' deduced from the completeness
of $J\left(\theta,r,A_{0},A_{1}\right)$, e.g., by using that completeness
to show that, for any sequence $\left\{ x_{m}\right\} _{m\in\mathbb{N}}$
in $\left\langle A_{0},A_{1}\right\rangle _{\theta,(r)}$ which satisfies
$\sum_{m=1}^{\infty}\left\Vert x_{m}\right\Vert _{\left\langle A_{0},A_{1}\right\rangle _{\theta,(r)}}<\infty$
there exists an element $y\in\left\langle A_{0},A_{1}\right\rangle _{\theta}$
such that $\lim_{N\to\infty}\left\Vert y-\sum_{m=1}^{N}x_{m}\right\Vert $$_{\left\langle A_{0},A_{1}\right\rangle _{\theta}}=0$.
\textcolor{blue}{}%

\bigskip{}

\textbf{\textit{Acknowledgment:}} It is a pleasure to express my warm
thanks to Svante Janson and to Mieczys\l{}aw Masty\l{}o for some very
helpful comments about the previous version of this note.


\begin{thebibliography}{10}
\bibitem{AvniE2014}E.~Avni, The periodic complex method in interpolation
spaces. \textit{Rev.~Mat.~Complut.} \textbf{27} (2014), 156--166.
Earlier version: \texttt{arXiv:1203.4183 {[}math.FA{]}}

\bibitem{BerghJLofstromJ1976} J.~Bergh and J.~Löfström,\textit{
Interpolation spaces. An Introduction.} Grundlehren der mathematische
Wissenschaften 223, Springer, Berlin-Heidelberg-New York, 1976.

\bibitem{CalderonA1964} A.~P.~Calderón, Intermediate spaces and
interpolation, the complex method. \textit{Studia Math.,} \textbf{24}
(1964), 113--190. \texttt{\href{http://matwbn.icm.edu.pl/ksiazki/sm/sm24/sm24110.pdf}{http://matwbn.icm.edu.pl/ksiazki/sm/sm24/sm24110.pdf}}

\bibitem{CwikelM1978-complex} M.\ Cwikel, Complex interpolation,
a discrete definition and reiteration, \textit{Indiana Univ.\ Math.\ J.}
\textbf{27} (1978), 1005--1009.

\bibitem{CwikelMJansonS2007} M.\ Cwikel and S.\ Janson, Complex
interpolation of compact operators mapping into the couple $(FL^{\infty},FL_{1}^{\infty})$.
In Contemporary Mathematics, \textbf{445}, (L.~De Carli and M.~Milman,\texttt{
}eds.). American Mathematical Society, Providence R.I., 2007, 71--92.
\textit{Preliminary version:} \texttt{arXiv:math/0606551v1 {[}math.FA{]}}.

\bibitem{CwKa}M.~Cwikel and N.~J.~Kalton, Interpolation of compact
operators by the methods of Calder\'on and Gustavsson\textendash{}Peetre,
\textit{Proc.~Edinburgh Math.~Soc.} \textbf{38} (1995), 261--276.
Earlier version: \texttt{arXiv:math.FA/9210206}

\bibitem{CwikelMKaltonNkmr} M.~Cwikel, N.~Kalton, M.~Milman and
R.~Rochberg, A Unified Theory of Commutator Estimates for a Class
of Interpolation Methods, \textit{Adv.~in Math. }\textbf{169} (2002),
241--312.

\bibitem{CwikelMMilmanMRochbergR2014B}M.~Cwikel, M.~Milman and
R.~Rochberg, Nigel Kalton and the interpolation theory of commutators,
\texttt{arXiv:1405.5686 {[}math.FA{]}}

\bibitem{GustavssonJPeetreJ1977} J.\ Gustavsson and J.\ Peetre,
Interpolation of Orlicz spaces, \textit{Studia Math.}\ \textbf{60}
(1977), 33--59. \texttt{\href{http://matwbn.icm.edu.pl/ksiazki/sm/sm60/sm6014.pdf}{http://matwbn.icm.edu.pl/ksiazki/sm/sm60/sm6014.pdf}}

\bibitem{JansonS1981}S.~Janson, Minimal and maximal methods of interpolation.
\textit{J.~Functional Analysis}, \textbf{44} (1981), 50--73.

\bibitem{KatznelsonY1968} Y.\ Katznelson, An introduction to Harmonic
Analysis, Wiley, New York-London-Sydney-Toronto, 1968.

\bibitem{LionsJPeetreJ1964} J.-L.~Lions and J.~Peetre, Sur une
classe d'espaces d'interpolation. \textit{Inst.\ Hautes Etudes Sci.\ Publ.\ Math.,}
\textbf{19}, (1964), 5--68.

\bibitem{OvchinnikovV1984} V.~I.~Ovchinnikov, The method of orbits
in interpolation theory. \textit{Mathematical Reports,} \textbf{1},
(1984), no.~2, pp.~i--x and 349-515.

\bibitem{PeetreJ1971} J.\ Peetre, Sur l'utilization des suites inconditionellement
sommables dans la théorie des espaces d'interpolation.\ \textit{Rend.\ Sem.\ Mat.\ Univ.\ Padova}
\textbf{46} (1971), 173--190. \end{thebibliography}
\end{document}